\documentclass[10pt]{amsart}

\theoremstyle{definition}
\newtheorem{thm}{Theorem}[section]
\newtheorem{lem}[thm]{Lemma}
\newtheorem{prp}[thm]{Proposition}
\newtheorem{dfn}[thm]{Definition}
\newtheorem{cor}[thm]{Corollary}

\newtheorem{rmk}[thm]{Remark}

\newtheorem{exa}[thm]{Example}
\newtheorem{pbm}[thm]{Problem}

\newtheorem{qst}[thm]{Question}
\newenvironment{pff}{{\em Proof:}}{\QED}

\newcommand{\beq}{\begin{equation}}
\newcommand{\eeq}{\end{equation}}
\newcommand{\beqr}{\begin{eqnarray*}}
\newcommand{\eeqr}{\end{eqnarray*}}
\newcommand{\limi}[1]{\lim_{{#1} \to \infty}}
\newcommand{\bit}{\begin{itemize}}
\newcommand{\eit}{\end{itemize}}
\newcommand{\ts}[1]{{\textstyle{#1}}}

\newcommand{\af}{\alpha}
\newcommand{\bt}{\beta}

\newcommand{\dt}{\delta}
\newcommand{\ep}{\varepsilon}
\newcommand{\zt}{\zeta}
\newcommand{\et}{\eta}

\newcommand{\io}{\iota}
\newcommand{\te}{\theta}
\newcommand{\ld}{\lambda}

\newcommand{\ph}{\varphi}

\newcommand{\rh}{\rho}
\newcommand{\om}{\omega}
\newcommand{\ta}{\tau}

\newcommand{\Q}{{\mathbf{Q}}}
\newcommand{\Z}{{\mathbf{Z}}}
\newcommand{\R}{{\mathbf{R}}}
\newcommand{\C}{{\mathbf{C}}}
\newcommand{\N}{{\mathbf{N}}}

\newcommand{\id}{{\mathrm{id}}}

\newcommand{\diam}{{\mathrm{diam}}}

\newcommand{\spec}{{\mathrm{sp}}}

\newcommand{\andeqn}{\,\,\,\,\,\, {\mbox{and}} \,\,\,\,\,\,}
\newcommand{\A}{\andeqn}
\newcommand{\QED}{\rule{0.4em}{2ex}}

\newcommand{\csmf}{compact smooth manifold}
\newcommand{\ca}{C*-algebra}
\newcommand{\ct}{continuous}
\newcommand{\pj}{projection}
\newcommand{\mf}{manifold}

\newcommand{\nbhd}{neighborhood}

\newcommand{\ifo}{if and only if}

\newcommand{\hme}{homeomorphism}
\newcommand{\dfe}{diffeomorphism}
\newcommand{\mh}{minimal homeomorphism}
\newcommand{\md}{minimal diffeomorphism}
\newcommand{\tgca}{transformation group \ca}

\newcommand{\cms}{compact metric space}
\newcommand{\chs}{compact Hausdorff space}
\newcommand{\cfn}{continuous function}

\newcommand{\rsz}[1]{\raisebox{0ex}[0.8ex][0.8ex]{$#1$}}

\title[Crossed products by minimal diffeomorphisms]{When
are crossed products by minimal diffeomorphisms isomorphic?}

\author{N.\  Christopher Phillips}

\address{Department of Mathematics, University  of Oregon,
       Eugene OR 97403-1222, USA.}

\email[]{ncp@darkwing.uoregon.edu}

\date{24 April 2002}  

\thanks{
  \\
\indent
Research partially supported by NSF grants DMS 9706850
  and DMS 0070776, and by the MSRI}

\begin{document}

\setcounter{section}{-1}

\begin{abstract}
We discuss the isomorphism problem for both C* and smooth
crossed products by \md s.
For C* crossed products, examples demonstrate the failure of the obvious
analog of the Giordano-Putnam-Skau Theorem on \mh s of the Cantor set.
For smooth crossed products, there are many open problems.
\end{abstract}

\maketitle

\section{Introduction}

\indent
A remarkable theorem of Giordano, Putnam, and Skau~(\cite{GPS};
see Theorem~\ref{GPSThm} below)
gives a dynamical characterization
of isomorphism of the \tgca s of \mh s of the Cantor set.
The proof depends, among other things, on the fact~\cite{El1} that
simple direct limits of circle algebras with real rank zero and
with the same scaled ordered K-theory are necessarily isomorphic.
Recent progress in the classification of simple \ca s~\cite{Ln}
and the structure of crossed products~\cite{LP2} (see the
survey~\cite{LP1b}) has made it possible in some cases to prove the
isomorphism of crossed products by \md s with the same Elliott
invariants.
Some examples have been constructed~\cite{Ph9},
which we survey in this paper.
The examples show is that the analog of the condition of~\cite{GPS}
is far too strong to correspond to isomorphism of the \ca s.
There are no clear candidates for the correct condition, but
the examples rule out a number of possibilities.

For a \md\  satisfying an additional technical condition, it is
possible to construct a smooth crossed product, which is a
locally multiplicatively convex Fr\'{e}chet algebra.
The smooth irrational rotation algebra is a well known example
of this construction.
The smooth crossed product presumably preserves much more information
about the dynamics than the C* crossed product, although so far
essentially no theorems to this effect are known.
It turns out that very little is known about smooth crossed products.
In the second half of this survey, we discuss conditions for the
existence of smooth crossed products, the isomorphism problem,
and some of the other interesting open questions.

This paper has four sections.
In the first, which may serve as a more extended introduction,
we discuss what is known in the case of \mh s of the Cantor
set and of the circle, another low dimensional case where the
situation can be completely described.
We then discuss four previously known examples which suggest, but do not
conclusively demonstrate, that the general case is rather different.
In the next section, we present four examples of pairs of \md s
which do not satisfy the condition of~\cite{GPS} but for which it has
recently become possible to prove that the crossed product \ca s
are in fact isomorphic.
The methods used to distinguish the \hme s are different in each case.
We raise a few specific questions about these examples, and give a
brief discussion of the problem of relating isomorphism of the
crossed product \ca s to the dynamics.
The third section describes a general sufficient condition for the
existence of a smooth crossed product, and shows that the
diffeomorphisms in at least some of our examples satisfy this
condition.
The problem, then, is to find a dynamical characterization
of isomorphism of smooth crossed products by \md s.
Smooth crossed products are also natural examples for
Connes's noncommutative geometry.
In the last section, we discuss some more elementary questions
about smooth crossed products which are still open, mostly about
the analogs of stable and real rank.
These arose when thinking about the isomorphism question for
smooth crossed products, and realizing how little is in fact known.
We give a recent example of Schweitzer which shows that at least
one of these questions can't be answered using the general theory
of smooth subalgebras of \ca s.

We are grateful to the organizers of the conference on 
Operator Algebras and Mathematical Physics at Constan\c{t}a in
July 2001, at which the results discussed in Section~2 were
presented, to Anatole Katok,
Larry Schweitzer, and Christian Skau for helpful email
correspondence, and to Ian Putnam and Larry Schweitzer for permission
to present here their unpublished examples.

\section{Flip conjugacy, orbit equivalence, and crossed product
  C*-algebras}\label{Intro}

\indent
The theorem of Giordano, Putnam, and Skau is as follows:

\begin{thm}[Theorem~2.1 of \cite{GPS}]\label{GPSThm}
Let $X$ be the Cantor set, and let $h_1, \, h_2 \colon X \to X$
be \mh s.
Then $C^* (\Z, X, h_1) \cong C^* (\Z, X, h_2)$ \ifo\  $h_1$ and $h_2$
are strong orbit equivalent.
\end{thm}

The precise definition of strong orbit equivalence is given in
Definition~1.3 of \cite{GPS}.
Since it is slightly technical and will not be needed,
we do not reproduce it here.
Rather, we define several related conditions,
one stronger than strong orbit equivalence and one weaker.
The first is very close to the obvious notion of isomorphism of
\hme s.

\begin{dfn}\label{FC}
Let $X_1$ and $X_2$ be topological spaces,
and let $h_1 \colon X_1 \to X_1$ and $h_2 \colon X_2 \to X_2$
be \hme s.
We say that $h_1$ and $h_2$ are {\emph{conjugate}} if there is a
\hme\  $g \colon X_1 \to X_2$
such that $g \circ h_1 \circ g^{-1} = h_2$.
We say that $h_1$ and $h_2$ are {\emph{flip conjugate}} if $h_1$ is
conjugate to either $h_2$ or $h_2^{-1}$.
\end{dfn}

If $h_1$ and $h_2$ are flip conjugate and $X_1$ and $X_2$
are locally compact, it is immediate that
$C^* (\Z, X_1, h_1) \cong C^* (\Z, X_2, h_2)$.

\begin{dfn}\label{TOE}
Let $X_1$ and $X_2$ be topological spaces,
and let $h_1 \colon X_1 \to X_1$ and $h_2 \colon X_2 \to X_2$
be \hme s.
We say that $h_1$ and $h_2$ are {\emph{topologically orbit equivalent}}
if there is a \hme\  $g \colon X_1 \to X_2$ such that,
for all $x \in X_1$,
\[
g ( \{ h_1^n (x) \colon n \in \Z \} )
   = \{ h_2^n (g (x)) \colon n \in \Z \}.
\]
\end{dfn}

That is, $g$ is required to map the orbits of the \hme\  $h_1$
exactly to the orbits of the \hme\  $h_2$.
This definition is adapted from a similar definition in
measurable dynamics.

For our purposes, the important facts about strong orbit equivalence
are that flip conjugacy implies strong orbit equivalence and that
strong orbit equivalence implies topological orbit equivalence.
There are many examples of \mh s of the Cantor set
which are strong orbit equivalent but not flip conjugate.
For example, flip conjugacy preserves
topological entropy on compact spaces
(Theorems~7.2 and~7.3 of~\cite{Wl}),
but all entropies in $[0, \infty]$ occur in every
strong orbit equivalence class
(Theorem~6.1 of \cite{Sg1}, Theorem~7.1 of \cite{Sg2},
and Theorem~7.1 of \cite{Sg3}).
However, even topological orbit equivalence preserves the space of
invariant probability measures.
See the proof of (i) implies (iii) in Theorem~2.2 of~\cite{GPS},
at the beginning of Section~5 there.

In the most interesting higher dimensional cases, the distinction
between the equivalence relations disappears.
Two orbit equivalent \mh s of a connected \cms\  are necessarily
flip conjugate.
See Theorem~3.1 and Remark~3.4 of \cite{BT},
or Proposition~5.5 of~\cite{LP1}.
Nevertheless, for \mh s of the circle $S^1$, it is true that
isomorphism of the \tgca s implies flip conjugacy.
This follows from the fact that every \mh\  of $S^1$ is conjugate
to an irrational rotation
(Proposition~11.1.4 and Theorem 11.2.7(1) of~\cite{KH}),
and the computation of the scaled ordered
K-theory of the irrational rotation \ca s in the Appendix of~\cite{PV},
which shows that $K_0 (A_{\te_1}) \not\cong K_0 (A_{\te_2})$
unless $\te_1$ has the same image as $\pm \te_2$ in $\R / \Z$.

Evidence has accumulated that strong orbit equivalence
(or flip conjugacy) is not the relation on general \mh s
which corresponds to isomorphism of the \tgca s.
We describe four known suggestive examples.
(We should also note that flip conjugacy has for some time been
known to fail in Examples~\ref{Rouhani} and~\ref{JiFurst} in the
next section.)
Since it plays a crucial role in the discussion, we give a formal
definition of the Elliott invariant of a unital \ca.
See, for example,~\cite{El2}.

\begin{dfn}\label{EllInv}
Let $A$ be a separable unital \ca.
Its {\emph{Elliott invariant}} consists of:
\begin{itemize}
\item
The abelian group $K_1 (A)$.
\item
The scaled ordered abelian group $K_0 (A)$, in which the scale is
the distinguished element $[1]$ and the order is defined by
$\et > 0$ \ifo\  there are an integer $n$ and a \pj\  $p \in M_n (A)$
such that $\et = [p]$.
\item
The simplex $T (A)$ of tracial states on $A$, equipped with the weak*
topology.
\item
The pairing $T (A) \times K_0 (A) \to \R$ determined by
$(\ta, \, [p]) \mapsto \ta (p)$.
\end{itemize}
An {\emph{isomorphism}} from the Elliott invariant of $A$ to that of
$B$ consists of a group isomorphism $\ph_1 \colon K_1 (A) \to K_1 (B)$,
an isomorphism
$\ph_0 \colon K_0 (A) \to K_0 (B)$ of scaled ordered groups,
and an isomorphism $f \colon T (B) \to T (A)$ of simplexes,
such that the pairs $(\ta, \, \ph_0 (\et))$ and $(f (\ta), \, \et)$ have
the same image in $\R$ for all $\ta \in T (B)$ and $\et \in K_0 (A)$.
\end{dfn}

The unital case  of the Elliott conjecture asserts that if
two simple separable nuclear (but not type I) \ca s have isomorphic
Elliott invariants, then the \ca s are isomorphic.
There is good evidence for this conjecture in ``low rank'' cases.
It holds for the purely infinite case under the single additional
assumption that the algebras satisfy the Universal Coefficient Theorem
\cite{Kr}, \cite{Ph}.
Among the many results in the stably finite case,
we cite~\cite{EGL} and~\cite{Ln}.

We now describe the examples.
The first is about 15 years old.

\begin{exa}\label{SmNCTorus}
In the discussion on Pages 506--507 of~\cite{BCEN}, it is shown
that there are two nonisomorphic antisymmetric bicharacters
$\rh$ and $\rh'$ on $\Z^3$ which are injective on $\Z^3 \wedge \Z^3$
and such that the Elliott invariants of the corresponding
higher dimensional noncommutative toruses $A_{\rh}$ and
$A_{\rh}'$ are isomorphic.
It follows from more recent work (\cite{LQ} and~\cite{El1})
that $A_{\rh} \cong A_{\rh}'$.
However, Theorem~2 of~\cite{BCEN} shows that the
corresponding smooth subalgebras are not isomorphic.
\end{exa}

This example, although suggestive,
does not directly bear on the question, since these algebras
are not obviously crossed products.
Even if they were,
nonisomorphism of smooth subalgebras would only obviously rule out
smooth flip conjugacy.

In preparation for the next example,
recall that an {\emph{extension}} of a \hme\  $h \colon X \to X$
of a compact metric space $X$ consists of a compact metric space $Y$,
a \hme\  $k \colon Y \to Y$, and a \ct\  surjective map
$\pi \colon Y \to X$ such that $\pi \circ k = h \circ \pi$, and further
recall (see, for example, page~157 of~\cite{As}) that the extension
is {\emph{almost one to one}} if there exists a point $x \in X$
such that $\pi^{-1} (x) \subset Y$ consists of just one point.
(If $k$ is minimal, then in fact $\pi^{-1} (x)$ will be a one point
set for ``most'' $x \in X$.)

\begin{exa}\label{GjJo}
It is shown in Theorem~4 of~\cite{GJ} that every \mh\  $h$ of
the Cantor set $X$ has an almost one to one extension which is a
\mh\  $k$ of a nonhomogeneous space $Y$ such that $C^* (\Z, X, h)$ and
$C^* (\Z, Y, k)$ have isomorphic Elliott invariants.
Since it is not homogeneous, the space $Y$ is not
homeomorphic to $X$, and in particular $h$ and $k$ can't possibly
be strong orbit equivalent, or even topologically orbit equivalent.

The space $Y$ has covering dimension $\dim (X) = 1$.
Theorem~\ref{Isom} below does not apply, because $X$ and $Y$ are
not manifolds.
However, it is probably now easy to prove in this case that isomorphism
of the Elliott invariants implies isomorphism of the \ca s.
See Remark~\ref{NonMfIso} below.
\end{exa}

One might argue that one should only consider \mh s of the same
space, or at least of spaces of the same dimension.
However, almost one to one extensions are a standard
construction of closely related \hme s in dynamics.

Next, we present an unpublished example of Ian Putnam.
It is reproduced here with his permission.

\begin{exa}[Putnam]\label{Putnam}
For any $\te \in \R \setminus \Q$ let $r_{\te} \colon S^1 \to S^1$
be rotation by $2 \pi \te$.
Further let $g_{\te}$ be a \mh\  of a Cantor set $X_{\te} \subset S^1$
obtained from a Denjoy \hme\  $g_{\te}^{(0)} \colon S^1 \to S^1$
as follows~\cite{PSS}.
Any Denjoy \hme\  has a unique minimal set $X$,
which is homeomorphic to the Cantor set.
Choose $g_{\te}^{(0)}$ to have rotation number $\te$ and such that
the unique minimal set $X_{\te} \subset S^1$ has the property that
the image of $S^1 \setminus X_{\te}$ under the semiconjugation to
$r_{\te}$ is a single orbit of $r_{\te}$.
Let $g_{\te} = g_{\te}^{(0)} |_{X_{\te}}$.
See Section~3 of~\cite{PSS} for details, particularly Corollary~3.2
and Definitions~3.3 and~3.5, noting that we are requiring the set
$Q$ there to consist of exactly one orbit.

Following Remark~1 in Section~3 of \cite{PSS}, we may
consider the C* subalgebra of the bounded Borel functions on $S^1$
generated by $C (S^1)$ and the characteristic functions of all
sets $\exp ( 2 \pi i [n \te, (n + 1) \te ) )$ for $n \in \Z$.
This \ca\  has an automorphism $\bt_{\te}$ given by
rotation by $2 \pi \te$,
and we can take $X_{\te}$ to be its maximal ideal space and
$g_{\te}$ to be the \hme\  determined by $\bt_{\te}$.

Now let $\te_1, \, \te_2 \in \R \setminus \Q$ be numbers such that
$1, \, \te_1, \, \te_2$ are rationally independent.
Define \hme s
\[
h_1 = r_{\te_1} \times g_{\te_2}
    \colon S^1 \times X_{\te_2} \to S^1 \times X_{\te_2}
\]
and
\[
h_2 = r_{\te_2} \times g_{\te_1}
    \colon S^1 \times X_{\te_1} \to S^1 \times X_{\te_1}.
\]
It follows from Proposition~\ref{PtEllInv} below
that the crossed products by
these \hme s are simple \ca s with unique traces, and have
isomorphic Elliott invariants.

We show that the \hme s $h_1$ and $h_2$ are not topologically orbit
equivalent.
Suppose we had a topological orbit equivalence $f$.
By the proof of (i) implies (iii) in Theorem~2.2 of~\cite{GPS},
at the beginning of Section~5 there, $f$ preserves the invariant measures.
(The proof works without the
restriction that the space be the Cantor set.)
So the sets of possible measures of compact open subsets are the same
for both systems.
For $S^1 \times X_{\te_2}$ this set contains $\te_2$ and is
contained in $\Z + \te_2 \Z$,
while for $S^1 \times X_{\te_1}$ this set contains $\te_1$ and is
contained in $\Z + \te_1 \Z$.
This is a contradiction.

As in Example~\ref{GjJo}, Remark~\ref{NonMfIso} suggests that it
should be easy to adapt known results to prove that the crossed
product \ca s are isomorphic.
\end{exa}

Unlike in Example~\ref{GjJo}, in this example the spaces are
homeomorphic.

The last example is from our earlier work with Qing Lin.

\begin{exa}\label{OddSpheres}
For a \md\  of a sphere $S^n$ with $n \geq 3$ odd, the Elliott invariant
of the \tgca s is known.
(See Section~5 of~\cite{LP1} and Example~4.6 of~\cite{Ph8}.
In fact, the calculation works for \mh s, using Corollary~VI.12
of~\cite{Ex} in place of Corollary~3 in Section~5 of~\cite{Cn0}.)
It depends only on the simplex of invariant Borel probability measures,
and in particular is independent of $n$, as long as $n \geq 3$, and
of other properties of the \dfe.
It follows from Theorem~3 of~\cite{FH} that every odd sphere
admits a uniquely ergodic \md, and from~\cite{Wn} that every odd sphere
of dimension at least~$3$
admits a \md\  with any given finite number of ergodic measures.
If the Elliott conjecture holds for the corresponding \tgca s, then
those \md s having a given finite number of ergodic measures all
give isomorphic \ca s independent of $n$, as long as $n \geq 3$.
Moreover, it is likely~\cite{Kt}, although it remains unproved, that
an odd sphere in fact admits many nonconjugate uniquely ergodic \md s.
\end{exa}

Currently known classification theorems are not
adequate to prove isomorphism of these \ca s from isomorphism of their
Elliott invariants.
The theorem presently available requires real rank zero, but
these \ca s have no nontrivial \pj s
(Corollary~3 in Section~5 of~\cite{Cn};
Corollary~12 in Section~6 of~\cite{Ex}).
 
In the next section, we describe several examples of \md s of
compact connected manifolds for which it is now possible to prove
that the \tgca s are isomorphic, while the \dfe s are not
flip conjugate and hence not even topologically orbit equivalent.
In several of our examples, the \mf s on which the \dfe s act
are identical, in another they are different but have the same
dimension, and in another they have different dimensions.
The variety of different examples gives a collection of properties
of minimal \dfe s which are not invariants of the \tgca s.

The following result is a special case of results of~\cite{LP2}
and~\cite{Ln},
and is what we use to establish isomorphisms of crossed products.
All crossed products covered by it have real rank zero, by~\cite{LP2}
and~\cite{Ph8b}.

\begin{thm}\label{Isom}
For $j = 1, \, 2$ let $M_j$ be a \csmf, and let
$h_j \colon M_j \to M_j$ be a \md.
Assume that the maps
\[
K_0 ( C^* (\Z, M_j, h_j) )
 \to {\mathrm{Aff}} ( T ( C^* (\Z, M_j, h_j) ) ),
\]
from the $K_0$-groups to the spaces of real valued affine \cfn s
on the trace spaces, have dense range.
Assume that the Elliott invariants (Definition~\ref{EllInv}) of
$C^* (\Z, M_1, h_1)$ and $C^* (\Z, M_2, h_2)$ are isomorphic.
Then $C^* (\Z, M_1, h_1) \cong C^* (\Z, M_2, h_2)$.
\end{thm}

In the theorem, density of the image of $K_0$ in the affine
function space is equivalent to the \tgca\  having real rank zero.
See~\cite{Ph8b}.

\begin{rmk}\label{NonMfIso}
The main, although not the only, use of smoothness in the proof
of Theorem~\ref{Isom} is for an exponential length bound in~\cite{LP2}.
However, when the covering dimension of the space is at most~$2$,
exponential length bounds are easy.
See Section~2 of~\cite{PR}.
So current methods can probably be easily
modified to prove that Theorem~\ref{Isom}
holds for \mh s of \cms s with covering dimension at most~$2$,
possibly under the restriction that there be at most countably many
ergodic measures.
\end{rmk}

In our examples, we computed the K-theory
using the Pimsner-Voiculescu exact sequence, Theorem~\ref{PVSeq} below.
In most of them, we computed the order on $K_0$ using
Exel's rotation numbers for automorphisms~\cite{Ex}.
A relatively easy case is carried out in Example~4.9 of~\cite{Ph8}.

We finish this section by giving the computation relevant to
Example~\ref{Putnam}, since it has not appeared elsewhere.
The order computation is somewhat different, since Exel's methods
are easily applied only to \hme s of connected spaces.
For use here, and for later reference,
we first state the Pimsner-Voiculescu exact sequence~\cite{PV}
for the special case
of a crossed product of a compact space by a \hme.

\begin{thm}[Pimsner-Voiculescu exact sequence]\label{PVSeq}
Let $X$ be a \chs, and let $h \colon X \to X$ be a \hme.
Then there is a natural six term exact sequence
\[
\begin{array}{ccccc}
K^0 ( X )  & \stackrel{\id - h^*}{\longrightarrow}
    & K^0 ( X )  & \longrightarrow & K_0 ( C^* (\Z, X, h) )   \\
{\mathrm{exp}} \uparrow \hspace*{2em}  & & & &
            \hspace*{1em} \downarrow \partial       \\
K_1 ( C^* (\Z, X, h) ) & \longleftarrow & K^1 ( X ) &
    \stackrel{\id - h^*}{\longleftarrow} &
     K^1 ( X )
\end{array}.
\]
The maps $K^i (X) \to K_i ( C^* (\Z, X, h) )$ are the maps on
K-theory induced by the inclusion $C (X) \to C^* (\Z, X, h)$.
\end{thm}

\begin{prp}\label{PtEllInv}
Let $\af, \, \bt \in \R \setminus \Q$ be numbers such that
$1, \, \af, \, \bt$ are rationally independent.
Let $r_{\af} \colon S^1 \to S^1$ and
$g_{\bt} \colon X_{\bt} \to X_{\bt}$ be as in Example~\ref{Putnam}.
Let
\[
h = r_{\af} \times g_{\bt}
    \colon S^1 \times X_{\bt} \to S^1 \times X_{\bt}.
\]
Then $h$ is minimal and uniquely ergodic, and the Elliott
invariant of the crossed product
$A = C^* (\Z, \, S^1 \times X_{\bt}, \, h)$
is given as follows: $K_1 (A) \cong \Z^3$, there is a unique
tracial state $\ta$, and $\ta_*$ induces an order isomorphism
$K_0 (A) \to \Z + \af \Z + \bt \Z \subset \R$ (with the order on
the range given by restriction from $\R$) such that $\ta_* ([1]) = 1$.
\end{prp}

\begin{pff}
The restriction to $X_{\bt}$ of the semiconjugation map of
Corollary~3.2 of~\cite{PSS} is a \ct\  surjective map
$f \colon X_{\bt} \to S^1$ such that
$f \circ g_{\bt} = r_{\bt} \circ f$, and
there is a countable subset $T \subset X_{\bt}$ such that
$f |_{X_{\bt} \setminus T}$ is injective.
Therefore
\[
\id_{S^1} \times f \colon S^1 \times X_{\bt} \to S^1 \times S^1
\]
is a \ct\  surjective map such that
\[
(\id_{S^1} \times f) \circ h
 = (r_{\af} \times r_{\bt}) \circ (\id_{S^1} \times f),
\]
which is injective on the dense set $S^1 \times (X_{\bt} \setminus T)$.
The \dfe\  $r_{\af} \times r_{\bt}$ is known to be minimal
(Proposition~1.4.1 of~\cite{KH})
and uniquely ergodic (Theorem~6.20 of~\cite{Wl}), with the ergodic
measure $\nu$ being the product of two copies of Lebesgue measure.

It follows that $h$ is minimal.
Indeed, if $Z \subset S^1 \times X_{\bt}$ is a nonempty closed
invariant subset, then $(\id_{S^1} \times f) (Z)$
is a nonempty closed subset of $S^1 \times S^1$ which is invariant under
$r_{\af} \times r_{\bt}$, whence
$(\id_{S^1} \times f) (Z) = S^1 \times S^1$.
Therefore $Z$ contains $S^1 \times (X_{\bt} \setminus T)$,
whence $Z = S^1 \times X_{\bt}$.

It also follows that $h$ is uniquely ergodic.
Indeed, the existence of an invariant Borel probability measure $\mu$
follows from the existence of such measures for
$r_{\af}$ and $g_{\bt}$ (see Section~3 of~\cite{PSS}),
or by general results.
Moreover, if $\mu_0$ is any other invariant Borel probability measure,
then the measure on $S^1 \times S^1$ given by
$E \mapsto \mu_0 ((\id_{S^1} \times f)^{-1} (E))$ is a Borel probability
measure which is invariant under $r_{\af} \times r_{\bt}$,
and hence equal to $\nu$.
The Lebesgue measure of $f (T)$ is zero because $f (T)$ is countable,
so $\mu_0 (S^1 \times T) = 0$.
For any Borel set $E \subset S^1 \times X_{\bt}$, we therefore get
\begin{align*}
\mu_0 (E) & = \mu_0 (E \cap [S^1 \times (X_{\bt} \setminus T)])
  = \nu ( (\id_{S^1} \times f)
        (E \cap [S^1 \times (X_{\bt} \setminus T)]) ) \\
& = \mu (E \cap [S^1 \times (X_{\bt} \setminus T)]) = \mu (E).
\end{align*}

We now know that $A$ is simple and has a unique trace, say $\ta$.
We use Theorem~\ref{PVSeq} to compute $K_* (A)$.
We have
\[
K^0 ( S^1 \times X_{\bt} ) \cong K^1 ( S^1 \times X_{\bt} )
  \cong K^0 ( X_{\bt} ),
\]
and we can identify
\[
\id - h^* \colon
   K^i ( S^1 \times X_{\bt} ) \to K^i ( S^1 \times X_{\bt} ),
\]
for both $i = 0$ and $i = 1$, with
\[
\id - g_{\bt}^* \colon  K^0 ( X_{\bt} ) \to  K^0 ( X_{\bt} ).
\]
The proof of Lemma~6.1 of~\cite{PSS} gives
\[
{\mathrm{Ker}} (\id - g_{\bt}^*) \cong \Z
\A {\mathrm{Coker}} (\id - g_{\bt}^*) \cong \Z^2.
\]
Therefore the sequence of Theorem~\ref{PVSeq} breaks up into two
short exact sequences
\[
0 \longrightarrow \Z^2 \longrightarrow K_i (A) \longrightarrow \Z
   \longrightarrow 0
\]
for $i = 0$ and $i = 1$.
Both sequences must split, giving $K_0 (A) \cong K_1 (A) \cong \Z^3$.

We next determine what the trace does on $K_0 (A)$.
Define $B_1 = C^* (\Z, S^1, r_{\af})$ and
$B_2 = C^* (\Z, X_{\bt}, g_{\bt})$.
These \ca s are simple and have unique traces, say $\ta_1$ and $\ta_2$.
The coordinate \pj s induce equivariant
maps $C (S^1) \to C (S^1 \times X_{\bt})$ and
$C (X_{\bt}) \to C (S^1 \times X_{\bt})$, and hence injective maps
$\ph_1 \colon B_1 \to A$ and $\ph_2 \colon B_2 \to A$.
Naturality in the Pimsner-Voiculescu exact sequence gives a
commutative diagram, in which the vertical maps are sums
(so that $s (\et_1, \et_2) = (\ph_1)_* (\et_1) + (\ph_2)_* (\et_2)$
etc.), as follows:
\[
\begin{array}{ccccc}
K^0 ( S^1 ) \oplus K^0 ( X_{\bt} )
  & \stackrel{(\io_1)_* \oplus (\io_2)_*}{\longrightarrow}
     & K_0 ( B_1) \oplus K_0 (B_2 )
     & \stackrel{\partial_1 \oplus \partial_2}{\longrightarrow}
     & K^1 ( S^1 ) \oplus K^1 ( X_{\bt} )   \\
s_0 \downarrow \hspace*{1.3em}  & &
\hspace*{1em} \downarrow s & &
             \hspace*{1.2em} \downarrow s_1       \\
K^0 ( S^1 \times X_{\bt} )
  & \stackrel{\io_*}{\longrightarrow}
     & K_0 ( A )
     & \stackrel{\partial}{\longrightarrow}
     & K^1 ( S^1 \times X_{\bt} )
\end{array}.
\]

We claim that $s$ is surjective.
Let $\et \in K_0 ( A )$.
The proof of Lemma~6.1 of~\cite{PSS} shows that the kernel of
$\id - g_{\bt}^* \colon K^0 ( X_{\bt} ) \to K^0 ( X_{\bt} )$ is
generated by $[1] \in K_0 (C ( X_{\bt} ))$.
Therefore the kernel of
$\id - h^*
  \colon K^1 ( S^1 \times X_{\bt} ) \to K^1 ( S^1 \times X_{\bt} )$
is generated by the class of the unitary $v (\zt, x) = \zt$.
It follows that $\partial (\et) = n [v]$ for some $n \in \Z$.
Moreover, with $u \in K^1 (S^1)$ given by $u (\zt) = \zt$, we get
$[v] = s_1 ([u], \, 0)$.
Since $\id - r_{\af}^* = 0$, there is $\mu \in K_0 (B_1)$ such that
$\partial_1 (\mu) = [u]$.
Now $\partial (\et - s (n \mu, 0)) = 0$ by commutativity, so
$\et - s (n \mu, 0)$ is in the image of $\io_*$.
It is easy to check that $s_0$ is surjective, and it follows from
commutativity of the left square that $\et - s (n \mu, 0)$
is in the image of $s$.
Therefore $\et$ is in the image of $s$.
This completes the proof of the surjectivity of $s$.

Uniqueness of the traces gives
\[
\ta_* \circ s = \ta_* \circ (\ph_1)_* + \ta_* \circ (\ph_2)_*
  = (\ta_1)_* + (\ta_2)_*.
\]
Since the range of $(\ta_1)_*$ is $\Z + \af \Z$ and, by Theorem~5.3
of~\cite{PSS}, the range of $(\ta_2)_*$ is $\Z + \bt \Z$, it
follows that the range of $\ta_*$ is exactly $\Z + \af \Z + \bt \Z$.
Since $K_0 (A) \cong \Z^3$, it follows that $\ta_*$ is an
isomorphism onto its image.

It remains only to show that
$\ta_* \colon K_0 (A) \to \Z + \af \Z + \bt \Z$ is an order isomorphism.
This follows from Theorem~4.5(1) of~\cite{Ph8}.
\end{pff}

\section{Examples}\label{Exs}

\indent
We describe here four examples of pairs of different \md s giving
isomorphic crossed products.
The \md s in the pairs are distinguished in a variety of ways:
the property of having topologically quasidiscrete spectrum,
acting on \mf s of different dimensions or on nonhomeomorphic \mf s of
the same dimension, and inducing automorphisms of singular cohomology
which are not conjugate.
Details of these examples will appear in~\cite{Ph9}.
We also describe an example which has not yet been proved to exist but
whose existence seems likely.
In this case, the \dfe s are distinguished by the behavior of
$\limi{n} d (h^n (x), \, h^n (y))$ for distinct points $x$ and $y$ in
the \mf.
At the end of the section, we list several obvious questions related to
our examples, and give a brief discussion of the problem of finding a
dynamical condition for isomorphism of the crossed product \ca s.

\begin{exa}[Furstenberg transformations on $( S^1 )^2$]\label{Rouhani}
Rouhani has in \cite{Rh} exhibited a Furstenberg transformation on
the $2$-torus $S^1 \times S^1$ which does not have
topologically quasidiscrete spectrum.
(A \hme\  $h \colon X \to X$ is said to have
topologically quasidiscrete spectrum if the linear map $C (X) \to C (X)$,
given by $f \mapsto f \circ h$, has sufficiently many
``quasieigenfunctions'', a kind of generalized eigenvector.
See Section~1 of~\cite{Rh}.)

Let $\te \in [0, 1] \setminus \Q$ be an irrational number,
to be chosen below, and let $r \colon S^1 \to \R$ be a smooth function,
also to be chosen below.
Define $h_1, \, h_2 \colon S^1 \times S^1 \to S^1 \times S^1$ by
\[
h_1 (\zt_1, \zt_2)
 = \ts{ \left( e^{2 \pi i \te} \zt_1, \, \zt_1 \zt_2 \right)  }
\andeqn
h_2 (\zt_1, \zt_2)
 = \ts{ \left( e^{2 \pi i \te} \zt_1,
   \, e^{2 \pi i r (\zt_1) } \zt_1 \zt_2 \right) }
\]
for $(\zt_1, \zt_2) \in S^1 \times S^1$.
(The only difference is the extra factor $\exp (2 \pi i r (\zt_1) )$
in the definition of $h_2$.)

It is observed in~\cite{Rh} that the
affine Furstenberg transformation $h_1$ is
always minimal and uniquely ergodic and always
has topologically quasidiscrete spectrum, and it is also
shown how to choose $\te$ and $r$ so that $h_2$ is minimal and
uniquely ergodic but does not have topologically quasidiscrete spectrum.
It is only proved in~\cite{Rh} that $r$ is \ct, but in fact
the choices made there give a smooth function $r$,
so that $h_1$ and $h_2$ are both \dfe s.
The Elliott invariants of the \tgca s of this type are computed
in~\cite{Kd}; a much faster calculation using more machinery is
given in Example~4.9 of~\cite{Ph8}.
They turn out to
depend only on $\te$ and the space of invariant measures.
Moreover, the dense range hypothesis in Theorem~\ref{Isom} is
satisfied, and it follows that the \tgca s are isomorphic.

However, the property of having topologically quasidiscrete spectrum
is preserved by flip conjugacy.
So $h_1$ and $h_2$ are not flip conjugate, hence not
topologically orbit equivalent.
\end{exa}

\begin{exa}[Affine Furstenberg transformations
   on $( S^1 )^3$]\label{JiFurst}
In Section~6.1 of the unpublished thesis of R.\  Ji \cite{Ji},
two affine Furstenberg transformations on $( S^1 )^3$ were given
which have the same Elliott invariant but are not flip conjugate,
and the question was raised whether they have isomorphic \tgca s.
Ji was not able to compute the order on $K_0$; he only computed
the map on $K_0$ determined by the (unique) trace.
However, Theorem~4.5(1) of \cite{Ph8} implies
that the order on $K_0$ is that determined by the trace.
The calculation is more complicated than for Example~\ref{Rouhani},
because of the presence of torsion and more Bott elements.

Fix $\te \in [0, 1] \setminus \Q$ and $m, \, n \in \Z$ with
$0 < m < n$.
Then the two affine Furstenberg transformations on $( S^1 )^3$, given by
\[
(\zt_1, \zt_2, \zt_3) \mapsto
  \left( \exp (2 \pi i \te) \zt_1, \, \zt_1^m \zt_2,
              \, \zt_2^n \zt_3 \right)
\]
and
\[
(\zt_1, \zt_2, \zt_3) \mapsto
  \left( \exp (2 \pi i \te) \zt_1, \, \zt_1^n \zt_2,
              \, \zt_2^m \zt_3 \right)
\]
(the difference is that $m$ and $n$ have been exchanged),
are not topologically orbit equivalent but have isomorphic
crossed product \ca s.

The isomorphism of the \ca s is obtained from Theorem~\ref{Isom}
and the calculation of the ordered K-theory,
because any affine Furstenberg transformation is minimal and uniquely
ergodic.
If the first \dfe\  above is called $h$, then for any nonzero values
of $m$ and $n$ it turns out that
$K_0 (C^* (\Z, M, h))$ and $K_1 (C^* (\Z, M, h))$
are both isomorphic to $\Z^4 \oplus \Z / m \Z \oplus \Z / n \Z$,
and that the isomorphism of $K_0 (C^* (\Z, M, h))$ with this group can
be chosen in such a way that the unique trace $\ta$ induces
the map
\[
\ta_* ( r_1, r_2, r_3, r_4, s_1, s_2) = r_1 + \te r_3
\]
and $K_0 (C^* (\Z, M, h))_+$ is identified with
\[
\{ ( r_1, r_2, r_3, r_4, s_1, s_2)
    \in \Z^4 \oplus \Z / m \Z \oplus \Z / n \Z \colon
     r_1 + r_3 \te > 0 \}
   \cup \{ 0 \}.
\]

There is in \cite{Ji} no indication of the proof that these two
\dfe s are not flip conjugate.
This, however, can be obtained by examining their effect on
singular cohomology with integer coefficients.
We have $H^1 ( (S^1)^3; \, \Z) \cong \Z^3$, and
with respect to suitable bases the two \dfe s induce maps with matrices
\[
\left( \begin{array}{ccc}
  1 & m & 0 \\ 0 & 1 & n \\ 0 & 0 & 1  \end{array} \right)
 \andeqn
\left( \begin{array}{ccc}
  1 & n & 0 \\ 0 & 1 & m \\ 0 & 0 & 1  \end{array} \right).
\]
Over $\Z$, the first of these matrices is similar to neither the
second nor its inverse, which rules out flip conjugacy.
(However, the two matrices {\emph{are}} similar over $\Q$.)

If we further choose $m$ and $n$ to be relatively prime and with
$| m |, \, | n | \geq 2$, then
one can exhibit yet a third affine Furstenberg transformation
on $(S^1)^3$ which gives the same \tgca\  yet is not flip conjugate
to either of the first two, namely
\[
(\zt_1, \zt_2, \zt_3) \mapsto
  \left( \exp (2 \pi i \te) \zt_1, \, \zt_1^{m n} \zt_2,
              \, \zt_2 \zt_3 \right).
\]
The point is that $\Z / m \Z \oplus \Z / n \Z \cong \Z / m n \Z$,
but again the actions on $H^1 ( (S^1)^3; \, \Z)$ rule out
flip conjugacy.
\end{exa}

\begin{exa}[Minimal diffeomorphisms
   on $S^2 \times S^1$ and $( S^1 )^3$]\label{3DNonHomeo}
Let $M_1 = S^2 \times S^1$, and let $u \in U (C (M_1))$
be given by $u (x, \zt) = \zt$.
Adapting methods of \cite{FH}, we can prove that there exists a 
uniquely ergodic minimal diffeomorphism $h \colon M_1 \to M_1$, with
unique invariant Borel probability measure $\mu$,
which is homotopic to the identity map and such
that the rotation number of $[u]$ with respect to $h$ and $\mu$
(in the sense of~\cite{Ex}) has the form $\exp (2 \pi i \te)$
for some $\te \in [0, 1] \setminus \Q$.

Let $M_2 = ( S^1 )^3$, and define
$h_2 \colon M_2 \to M_2$ by
\[
h_2 (\zt_1, \zt_2, \zt_3)
 = \left( \exp (2 \pi i \te) \zt_1,
             \, \zt_1 \zt_2, \, \zt_2 \zt_3 \right)
\]
for $(\zt_1, \zt_2, \zt_3) \in ( S^1 )^3$.

The same methods as used for Example~\ref{JiFurst} enable one to
compute the Elliott invariant of $C^* (\Z, M_2, h_2)$.
A similar calculation, but with different algebraic topology,
computes the Elliott invariant of $C^* (\Z, M_1, h_1)$.
One gets
\[
K_0 ( C^* (\Z, M_1, h_1)) \cong K_0 ( C^* (\Z, M_2, h_2)) \cong \Z^4
\]
and
\[
K_1 ( C^* (\Z, M_1, h_1)) \cong K_1 ( C^* (\Z, M_2, h_2)) \cong \Z^4.
\]
Moreover, the orders on both $K_0$ groups turn out to be described as
follows: there are generators $\et_1$, $\et_2$, $\nu_1$, and $\nu_2$
such that the unique trace $\ta$ on the algebra satisfies
\[
\ta_* ( \et_1 ) = 1, \,\,\,\,\,\,
\ta_* ( \nu_1 ) = \te, \andeqn
\ta_* ( \et_2 ) = \ta_* ( \nu_2 ) = 0,
\]
and the positive cone of $K_0$ consists exactly of $0$ together
with the elements $\et$ such that $\ta_* (\et) > 0$.
In particular, the two Elliott invariants are isomorphic.

It follows from Theorem~\ref{Isom} that the two
crossed product \ca s are isomorphic.
However, the \dfe s can't possibly be topologically orbit equivalent
because the spaces on which they act are not homeomorphic.
\end{exa}

\begin{exa}[Minimal diffeomorphisms
   on manifolds of different dimensions]\label{DiffDims}
As discussed in Example~\ref{OddSpheres}, it is expected that the \tgca s
of \mh s of odd spheres of dimension at least $3$ depend up to
isomorphism only on the space of invariant probability measures.
Because these \ca s have no nontrivial \pj s, current machinery
is not able to prove this.
By forming the products of uniquely ergodic examples of this type
with a suitable irrational rotation on the circle,
we can produce examples to which current methods apply.

Use Theorem~3 (in Section~3.8) of \cite{FH} to find,
for each odd $n \geq 3$,
a uniquely ergodic \md\  $h_n^{(0)} \colon S^n \to S^n$.
The Lefschetz fixed point theorem (Theorem 4.7.7 of~\cite{Sp})
can be used to show that an orientation
reversing \dfe\  of an odd sphere must have a fixed point.
So our \dfe s are all orientation preserving and therefore homotopic
to the identity map.

Next, one proves using results from \cite{Pr} that if
$h \colon X \to X$ is a uniquely
ergodic \mh\  of a connected \cms\  $X$,
then there is a dense $G_{\dt}$-set $T \subset S^1$ such that, for every
$\ld \in T$, the \hme\  of $X \times S^1$ given by
$(x, \zt) \mapsto (h (x), \, \ld \zt)$ is minimal and uniquely ergodic.
Let $T_n$ be the dense $G_{\dt}$-set obtained for $h_n^{(0)}$,
and let $T$ be the intersection of these sets, which is still a
dense $G_{\dt}$-subset of $S^1$.
Choose $\te \in [0, 1] \setminus \Q$
such that $\exp (2 \pi i \te) \in T$.
For odd $n \geq 3$, define a uniquely ergodic
\md\  $h_n \colon S^n \times S^1 \to S^n \times S^1$ by
$h_n (x, \zt)
 = \left( \rsz{ h_n^{(0)} (x) }, \, \exp ( 2 \pi i \te) \zt \right)$.

Each $h_n$ is homotopic to the identity map, and K-theory does not
detect the difference between different odd spheres,
so the Pimsner-Voiculescu exact sequence gives the same K-groups
for all the crossed products $C^* (\Z, \, S^n \times S^1, \, h_n)$.
In particular, the $K_0$-groups are all $\Z^4$.
Using methods similar to those described in previous examples,
one shows that there is a set
of four generators of $K_0 ( C^* (\Z, \, S^n \times S^1, \, h_n) )$
whose traces are $1$, $\te$, $0$, and $0$.
As before, Exel's rotation numbers are used, and the fact that
$n \geq 3$ is used to show that there is essentially
only one source for a noninteger trace, namely a class in
$K_0 ( C^* (\Z, \, S^n \times S^1, \, h_n) )$ whose image in
$K_1 ( C (S^n \times S^1))$ is the class of the unitary
$(x, \zt) \mapsto \zt$.

It follows from Theorem~\ref{Isom} that the \ca s
$C^* (\Z, \, S^n \times S^1, \, h_n)$ are all isomorphic.
No two of these \dfe s can be topologically orbit equivalent,
since they act on \mf s of different dimensions.
\end{exa}

To complete our collection of examples, we give a brief description
of an example whose existence seems plausible but has not been
proved.

\begin{exa}[Extensions of Furstenberg transformations]\label{FurstExt}
Let $\te \in \R \setminus \Q$, and let
$h_1 \colon S^1 \times S^1 \to S^1 \times S^1$ be given by
\[
h_1 (\zt_1, \zt_2)
 = \ts{ \left( e^{2 \pi i \te} \zt_1, \, \zt_1 \zt_2 \right)},
\]
as in Example~\ref{Rouhani}.
Fix a point $z_0 \in S^1 \times S^1$.
We believe it should be possible to construct a surjective map
$f \colon S^1 \times S^1 \to S^1 \times S^1$ and a \md\  (or at
least a \mh) $h_2 \colon S^1 \times S^1 \to S^1 \times S^1$
with the following properties:
\begin{itemize}
\item
$h_1 \circ f = f \circ h_2$.
(Thus, $h_2$ is an extension of $h_1$.)
\item
If we let $T$ denote the orbit of $z_0$ under $h_1$, then
$f^{-1} (S^1 \times S^1 \setminus T)$ is dense in $S^1 \times S^1$,
and the restriction of $f$ to this set is injective.
(In particular, $h_2$ is an almost one to one extension of $h_1$,
as defined before Example~\ref{GjJo}.)
\item
For $n \in \Z$, the set $I_n = f^{-1} ( \{ h_1^n (z_0) \} )$
is homeomorphic to $[0, 1]$.
\item
$\limi{n} \diam (I_n) = 0$ and $\lim_{n \to - \infty} \diam (I_n) = 0$.
\item
$f$ commutes with the projection to the first coordinate.
\item
$f$ is a homotopy equivalence.
\end{itemize}

To construct $h_2$ and $f$, we replace each point in the orbit
of $z_0$ under $h_1$ by a ``vertical'' interval, starting with $z_0$,
then $h_1 (z_0)$ and $h_1^{-1} (z_0)$, etc., with the lengths of the
inserted intervals chosen to go to zero fast enough that the
resulting space is still a torus.
To see the possibility of replacing one point by an interval in a
continuous manner,
consider the map $f_0 \colon \R^2 \to \R^2$ given by
\[
f_0 (x, y) = \left\{ \begin{array}{ll}
     (x, \, y)   & \hspace{3em}  | x | \geq 1 \\
     (x, \, |x| y)   & \hspace{3em}
            {\mbox{$| x | \leq 1$ and $| y | \leq 1$}}  \\
     (x, \, |x| + 2 (y - 1) )  & \hspace{3em}
            {\mbox{$| x | \leq 1$ and $1 \leq y \leq 2 - |x|$}}  \\
     (x, \, - |x| + 2 (y + 1) )  & \hspace{3em}
            {\mbox{$| x | \leq 1$ and $-1 \geq y \geq -(2 - |x|)$}}  \\
     (x, \, y)   & \hspace{3em}
            {\mbox{$| x | \leq 1$ and $| y | \geq 2 - |x|$}}
    \end{array} \right..
\]
This map is the identity outside a compact set, injective off
$\{ 0 \} \times [-1, \, 1]$, and sends $\{ 0 \} \times [-1, \, 1]$
to $(0, 0)$.
Smooth versions exist, but are more complicated to describe.
Then $h_2$ must be defined to carry $I_n$ homeomorphically to
$I_{n + 1}$ for each $n$.
The details, especially for the smooth case, seem rather difficult to
carry out, but we conjecture that this can be done.

If $f$ and $h_2$ exist, then the fact that $f$ is a homotopy equivalence
can be used to show that the corresponding map from
$C^* (\Z, \, S^1 \times S^1, \, h_1)$ to
$C^* (\Z, \, S^1 \times S^1, \, h_2)$ induces an isomorphism of
Elliott invariants.
We show that $h_1$ and $h_2$ are not flip conjugate.
Let $d$ be the metric on $S^1 \times S^1$ given by the maximum
of the differences of the two coordinates.
If $z_1, \, z_2 \in S^1 \times S^1$ are distinct, then
\[
\liminf_{n \to \infty} d (h_1^n (z_1), \, h_1^n (z_2)) > 0
\A \liminf_{n \to - \infty} d (h_1^n (z_1), \, h_1^n (z_2)) > 0.
\]
(Consider separately the cases in which $z_1$ and $z_2$ have
different or equal first coordinates.)
That is, $h_1$ is distal as defined at the beginning of Chapter~5
of~\cite{As}.
However, if $z_1$ and $z_2$ are any two points in
$f^{-1} ( \{ z_0 \} )$, then
\[
\lim_{n \to \infty} d (h_2^n (z_1), \, h_2^n (z_2)) = 0
\A \lim_{n \to - \infty} d (h_2^n (z_1), \, h_2^n (z_2)) = 0.
\]
These properties are unchanged if $d$ is replaced by an equivalent
metric or one of the \hme s is replaced by its inverse.
Therefore $h_1$ and $h_2$ are not flip conjugate,
and so can't be topologically orbit equivalent either.
\end{exa}

These examples leave open a number of interesting problems.

\begin{pbm}\label{AllIrrat}
Make a more careful modification of the work of~\cite{FH}, so as
to be able to construct examples like Example~\ref{3DNonHomeo} and
Example~\ref{DiffDims} for arbitrary irrational values of $\te$.
\end{pbm}

\begin{pbm}\label{DiffOnOddSph}
Do there exist essentially different uniquely ergodic \md s of the same
odd sphere?
Can they be used to produce essentially different versions of
Example~\ref{DiffDims}?
\end{pbm}

Among the examples of \md s we know, those for which
the crossed product has real rank zero are all uniquely ergodic.

\begin{qst}\label{NonUniqErgQ}
Is there a \md\  $h$ of a connected \csmf\  $M$ which is not uniquely
ergodic but for which the map
\[
K_0 ( C^* (\Z, M, h) )
 \to {\mathrm{Aff}} ( T ( C^* (\Z, M, h) ) ),
\]
as in Theorem~\ref{Isom}, still has dense range?
\end{qst}

Following the methods of~\cite{Ex}, one needs among other things
a \cfn\  $u \colon M \to S^1$ which has different rotation numbers
with respect to different invariant Borel probability measures on $M$.

Finally, we turn to the question which motivated the construction of
these examples.

\begin{pbm}\label{WhenCStarIso}
What dynamical relation on \mh s of \cms s, or \md s of \csmf s,
corresponds to isomorphism of the \tgca s?
\end{pbm}

At the moment, we don't even have any plausible candidates.
One might think of considering flow equivalence, as
discussed, for example, at the beginning of Section~1 of \cite{Pc2} and
in Definition~1.1 and the following discussion in \cite{Pc1}.
(The right notion would actually be flip flow equivalence.)
We do not know whether the \md s of Example~\ref{Rouhani}
and of Example~\ref{JiFurst} are flow equivalent.
Those in Example~\ref{3DNonHomeo} and in Example~\ref{DiffDims}
are not, because Theorem~2 of~\cite{Sz} implies that if
$h_1$ and $h_1$ are flow equivalent \dfe s on manifolds $M_1$ and $M_2$
then the universal covers of $M_1$ and $M_2$ are homeomorphic.
Moreover, flow equivalence of \md s only implies stable isomorphism:
see Section~2 of~\cite{Pc2}.

\begin{pbm}\label{FlowEqQ}
Are the \md s of Example~\ref{Rouhani} flow equivalent?
Are those of Example~\ref{JiFurst} are flow equivalent?
\end{pbm}

The following question arose in discussions of possible answers
to Problem~\ref{WhenCStarIso}.

\begin{pbm}\label{CommonExt}
For $j = 1, \, 2$ let $X_j$ be a \cms, and let
$h_j \colon X_j \to X_j$ be a \mh.
Suppose that $C^* (\Z, X_1, h_1) \cong C^* (\Z, X_2, h_2)$.
Does it follow that $h_1$ and $h_2$ have a common minimal almost
one to one extension (as defined before Example~\ref{GjJo})?
\end{pbm}

This is true in Example~\ref{GjJo}, since $k$ is already an
almost one to one extension of $h$.
It is true in Example~\ref{Putnam}, since $g_{\te_1} \times g_{\te_2}$
is a common minimal almost one to one extension.

This condition certainly does not imply isomorphism of the \tgca s.
In the notation of Example~\ref{Putnam}, the \hme\  $g_{\te}$
is a minimal almost one to one extension of $r_{\te}$.
However,
\[
K_1 ( C^* ( \Z, S^1, r_{\te} )) \cong \Z^2
\A
K_1 ( C^* ( \Z, X_{\te}, g_{\te} )) \cong \Z.
\]

If one wants extensions to preserve K-theory of the crossed products,
then some restriction on the space of the extension is necessary.
The following seems like a good test case.

\begin{pbm}\label{RestrCommonExt}
Let $h_1, \, h_2 \colon S^1 \times S^1 \to S^1 \times S^1$ be
Furstenberg transformations as in Example~\ref{Rouhani},
with $\te \in \R \setminus \Q$ arbitrary and
$r \colon S^1 \to \R$ an arbitrary smooth function.
Does there exist a common extension of
$h_1$ and $h_2$ which is a \mh\  of $S^1 \times S^1$?
\end{pbm}

\section{Smooth crossed products}

The examples in Section~\ref{Exs} show that \ca\  crossed products
preserve little information about \mh s of connected \cms s.
For the case of a diffeomorphism satisfying an additional condition,
one can construct instead a smooth crossed product.
It is natural to hope, especially in view of Example~\ref{SmNCTorus},
that the smooth crossed product might preserve more information.
In fact, very little is known about smooth crossed products by \md s;
even some very basic questions are open.
In this section and the next, we discuss the smooth crossed products
and raise some of these questions.

The foundations of the abstract theory of smooth crossed products of
Banach and Fr\'{e}chet algebras are laid in \cite{Sw2}.
Here, we are looking at smooth crossed products by actions of
$\Z$ on the Fr\'{e}chet algebra $C^{\infty} (M)$ of smooth functions on
a compact smooth manifold $M$, with the topology of uniform convergence
of all derivatives.
We can give a collection of seminorms on $C^{\infty} (M)$
which determines the topology as follows.
Choose finitely many smooth vector fields $X_1, X_2, \dots, X_N$ on
$M$ such that, at every $x \in M$, the vectors
$X_1 (x), \, X_2 (x), \, \dots, \, X_N (x)$
span the tangent space $T_x M$.
For $f \in C^{\infty} (M)$ and $n \geq 0$, define
\[
\| f \|_n
 = \sum_{1 \leq k_1, \dots, k_n \leq N}
    \| X_{k_n} X_{k_{n - 1}} \cdots X_{k_1} f \|_{\infty}.
\]
These are not submultiplicative, but using the product rule for
derivatives, one can show that there are enough submultiplicative
combinations of them to define the topology on $C^{\infty} (M)$.
The well behaved smooth crossed product by $\Z$ is the space of
${\mathcal{S}} (\Z, \, C^{\infty} (M), \, h)$
of sequences $s \colon \Z \to C^{\infty} (M)$ which decay rapidly
at infinity in each of the seminorms on $C^{\infty} (M)$.
Its topology is determined by the seminorms
\[
\| s \|_{n, d} = \sum_{k \in \Z} (1 + | k |)^d \| s (k) \|_n
\]
for $d, n \geq 0$.
Unfortunately, if $h$ is an arbitrary \dfe\  of $M$, then
${\mathcal{S}} (\Z, \, C^{\infty} (M), \, h)$ need not be an algebra
under convolution.
We will require that the action of $h$ be m-tempered in the
sense of Definition~3.1.1 of~\cite{Sw2} (see below),
which by Theorem~3.1.7 of~\cite{Sw2} implies that
${\mathcal{S}} (\Z, \, C^{\infty} (M), \, h)$ is in fact a
locally multiplicatively convex Fr\'{e}chet algebra.
In the case at hand, it will in fact be a *-algebra.
The m-temperedness condition is stated in terms of submultiplicative
seminorms on $C^{\infty} (M)$, but it is sufficient to require it for
the seminorms
$\| \cdot \|_0 + \| \cdot \|_1 + \cdots + \| \cdot \|_n$ with the
seminorms $\| \cdot \|_m$ on $C^{\infty} (M)$ being as introduced above.
This leads us to the following definition.

\begin{dfn}\label{Tempered}
Let $M$ be a compact smooth manifold, and let $h \colon M \to M$
be a \dfe.
We say that $h$ is {\emph{tempered}} if for every $m$ there is
$C > 0$ and $r \in \N$ such that for every $f \in C^{\infty} (M)$
and $n \in \Z$, we have
\[
\| f \circ h^n \|_m
  \leq C (1 + | n |)^r ( \| f \|_0 + \| f \|_1 + \cdots + \| f \|_m).
\]
\end{dfn}

This definition is really a condition on the
derivatives of $h^n$ as $| n | \to \infty$: they should grow
at most polynomially in $n$.
The best estimate for general \dfe s allows exponential growth
of the derivatives of $h^n$; see Example~\ref{NonT} below.
To formulate our condition in terms of the derivatives of $h^n$
requires setting up more notation than we want to introduce here.
In all the explicit examples we actually discuss, $M$ will be
$(S^1)^d \cong (\R / \Z)^d$ for some $d$,
and in this case, as we now explain, we can describe the situation
in terms of ordinary derivatives on $\R^d$.

Let $M = (\R / \Z)^d$, and let $h \colon M \to M$ be a smooth function.
Then $h$ has a universal cover ${\widetilde{h}} \colon \R^d \to \R^d$.
It is not unique, but
every other choice has the form $x \mapsto {\widetilde{h}} (x) + l$
for some $l \in \Z^d$.
Moreover, for every $k \in \Z^d$ there is $l \in \Z^d$
such that ${\widetilde{h}} (x + k) = {\widetilde{h}} (x) + l$ for
all $x \in \R^d$.
(This is just continuity and the condition that ${\widetilde{h}}$ gives
a well defined function $(\R / \Z)^d \to (\R / \Z)^d$.)
So the partial derivatives of ${\widetilde{h}}$,
of order at least~$1$, depend only on $h$ and
are periodic with period~$1$ in each coordinate.
In particular, they are bounded.
Letting ${\widetilde{h}}_i \colon \R^d \to \R$ be the $i$-th
component of ${\widetilde{h}}$,
and letting $D_k$ denote partial differentiation on $\R^d$ with respect
to the $k$-th coordinate,
we can therefore define
\[
\rh_m (h)
  = \max_{1 \leq i \leq d} \sum_{j_1 + \cdots + j_r = m}
   \| D_1^{j_1} D_2^{j_2} \cdots D_r^{j_r} {\widetilde{h}}_i \|_{\infty}
  \in [0, \infty).
\]

\begin{lem}\label{TCond}
Let $h \colon (\R / \Z)^d \to (\R / \Z)^d$ be a \dfe.
Let $\rh_m$ be as above.
Suppose that for every $m \geq 1$
there is $C > 0$ and $r \in \N$ such that for every $n \in \Z$ we have
$\rh_m (h^n) \leq C (1 + | n |)^r$.
Then $h$ is tempered.
\end{lem}

\begin{pff}
Let $M = (\R / \Z)^d$.
To every $f \in C^{\infty} (M)$ corresponds a function
${\widetilde{f}} \in C^{\infty} (\R^d)$
which is periodic with period~$1$ in each coordinate.
Moreover,
$(f \circ h){\widetilde{ }} = {\widetilde{f}} \circ {\widetilde{h}}$.
With an obvious choice of vector fields on $M$, we get
\[
\| f \|_m
 = \sum_{j_1 + \cdots + j_r = m}
    \| D_1^{j_1} D_2^{j_2} \cdots D_r^{j_r} {\widetilde{f}} \|_{\infty}.
\]

We now consider derivatives of ${\widetilde{f}} \circ {\widetilde{h}}$.
Applying the chain rule and the product rule, we find that a
derivative
\[
D_1^{j_1} D_2^{j_2} \cdots D_r^{j_r}
    ({\widetilde{f}} \circ {\widetilde{h}} ) (x),
\]
with $j_1 + j_2 + \cdots + j_r = m$, is a finite sum of products
of at most $m + 1$ terms, one of which is a partial derivative of
$f$ of order at most $m$ evaluated at ${\widetilde{h}} (x)$
and the rest of which are partial derivatives of components of
${\widetilde{h}}$ of order between $1$ and $m$ evaluated at $x$.
In particular, there is a constant $C_{m, d}$, not depending on
$f \in C^{\infty} (M)$ or $h \colon M \to M$, such that
\[
\| f \circ h \|_m
  \leq C_{m, d} ( \| f \|_0 + \| f \|_1 + \cdots + \| f \|_m)
    [1 + \rh_1 (h) + \rh_2 (h) + \cdots + \rh_m (h)]^m.
\]
Applying this to $h^n$ in place of $h$,
and using the estimate on $\rh_m (h^n)$ in the hypotheses, we get
a bound of the required sort for $n \geq 0$.
For $n \leq 0$, apply the same argument to $h^{-1}$.
\end{pff}

\begin{rmk}\label{OftenT}
For $h \colon \R / \Z \to \R / \Z$ and $n \geq 1$, we have
\[
({\widetilde{h}}^n)' (t)
 = {\widetilde{h}}' ({\widetilde{h}}^{n - 1} (t)) \cdot
   {\widetilde{h}}' ({\widetilde{h}}^{n - 2} (t))
   \cdots {\widetilde{h}}' (h (t)) \cdot {\widetilde{h}}' (t).
\]
The naive estimate therefore gives $\rh_1 (h^n) \leq \rh_1 (h)^n$.
If $h$ is a \dfe, then it is possible to have
$| {\widetilde{h}}' (t) | \geq 1$ everywhere only if
$| {\widetilde{h}}' (t) | = 1$ everywhere.
If furthermore $h$ is minimal, then one can hope that the
iterates $x, \, h (x), \, h^2 (x), \, \dots$ are distributed
well enough between places where the derivative is small and
where it is large that the overall growth of the derivative is
less than exponential.
In fact, one expects \cite{Kt} that it is reasonably
common for \md s of compact \mf s to be tempered.
\end{rmk}

\begin{exa}\label{IrratRot}
We show that rotations are tempered.
Fix $\te_1, \te_2, \dots, \te_d \in \R$.
Define $h \colon (S^1)^d \to (S^1)^d$ by
\[
h (\zt_1, \zt_2, \dots, \zt_d)
  = \left( e^{2 \pi i \te_1} \zt_1, \, e^{2 \pi i \te_2} \zt_2,
             \, \dots, \, e^{2 \pi i \te_d} \zt_d \right).
\]
Then we can take ${\widetilde{h}}$ to be
\[
{\widetilde{h}} (x_1, x_2, \dots, x_d)
  = (x_1 + \te_1, \, x_2 + \te_2, \, \dots, \, x_d + \te_d),
\]
so
\[
{\widetilde{h}}^n (x_1, x_2, \dots, x_d)
  = (x_1 + n \te_1, \, x_2 + n \te_2, \, \dots, \, x_d + n \te_d).
\]
It follows that $\rh_1 (h^n) = d$
and $\rh_m (h^n) = 0$ for all $m \geq 2$.
So $h$ is tempered.
\end{exa}

\begin{exa}\label{NonT}
Let $g \colon \R \to \R$ be a $C^{\infty}$ function with period $1$,
with $g' (t) > -1$ for all $t$, with $g (0) = 0$, and with $g' (0) = 1$.
Let $h \colon \R / \Z \to \R / \Z$ be the \dfe\  determined by
${\widetilde{h}} (t) = t + g (t)$.
Then for $n \geq 1$ we have $({\widetilde{h}}^n)' (0) = 2^n$.
So $h$ does not satisfy the conditions of Lemma~\ref{TCond},
and in fact it is easy to see that $h$ is not tempered.
\end{exa}

The \dfe\  of Example~\ref{NonT} is not minimal.
However, minimal examples are known to exist.
We use the minimal real analytic \dfe\  $F_{\af}$ with nonzero
topological entropy constructed in~\cite{Hr}.
For us, the relevant properties are in Proposition~5.1 and
Theorem~5.3 of~\cite{Hr}.
In particular, the proof of Proposition~5.1 of~\cite{Hr} depends on
the inequality, for a suitable invariant probability measure $\mu$,
\[
\inf_{n \in \N} \frac{1}{n} \int_M \log ( \| T F_{\af}^n (x) \| ) \,
  d \mu (x) > 0,
\]
in which $T F_{\af}^n$ is the tangent map of $F_{\af}^n$.
This inequality is incompatible with a polynomial bound on the
growth of the derivatives of powers of $F_{\af}$.

\begin{exa}\label{AffFurstT}
The affine Furstenberg transformations of Example~\ref{JiFurst}
are tempered.
To prove this, take
\[
{\widetilde{h}} (x_1, x_2, x_3)
  = (x_1 + \te, \, x_2 + m x_1, \, x_3 + n x_2).
\]
We can write this in the more general form
${\widetilde{h}} (x) = t + x + N x$,
where $t \in \R^d$ is fixed and $N \in M_d (\R)$ is nilpotent
with $N^d = 0$.
(Here $d = 3$.)
By induction, we find that if $n > 0$ then
\[
{\widetilde{h}}^n (x) = t + (1 + N)t + \cdots + (1 + N)^{n - 1} t
    + (1 + N)^n x.
\]
For suitable $t_n \in \R^d$, using $N^d = 0$, and with $\binom{n}{k}$
denoting the binomial coefficient,
\[
{\widetilde{h}}^n (x)
  = t_n + \sum_{k = 0}^{\min (n, \, d - 1) } \binom{n}{k} N^k x.
\]
All partial derivatives of all components
$( {\widetilde{h}}^n )_i$ of order~$1$ are constant,
and all partial derivatives of higher order are zero.
Since the binomial coefficients $\binom{n}{k}$, with $k \leq d - 1$,
are polynomials in $n$ of degree at most $d - 1$, it follows that
there is a constant $C$ such that
\[
| D_j ( {\widetilde{h}}^n )_i (x) | \leq C (1 + n)^{d - 1}
\]
for all $i$ and $j$, all $n > 0$, and all $x \in \R^d$.
The \dfe\  ${\widetilde{h}}^{-1}$ has the same form,
so the same method applies.
Therefore $h$ is tempered by Lemma~\ref{TCond}.
\end{exa}

The argument of Example~\ref{AffFurstT} actually applies to all of
the affine Furstenberg transformations of~\cite{Ji}.

\begin{exa}\label{RouhT}
The Furstenberg transformations of Example~\ref{Rouhani} are tempered.
For $h_1$, this is the same as Example~\ref{AffFurstT}.
For $h = h_2$, we take
\[
{\widetilde{h}} (x, y)
   = (x + \te, \, y + x + {\widetilde{r}} (x)),
\]
where ${\widetilde{r}} (t) = r (e^{2 \pi i t} )$.
By induction on $n$, we get
\[
{\widetilde{h}}^n (x, y)
  = \left( x + n \te, \,\, y + n x
           + {\textstyle{ \frac{1}{2} }} n (n - 1) \te
           + \sum_{k = 0}^{n - 1} {\widetilde{r}} (x + k \te) \right)
\]
for $n \geq 0$.
The first component has
$D_1 ( {\widetilde{h}}^n )_1 (x, y) = 1$,
while all other partial derivatives, of all orders, are zero.
The second component has
\[
D_2 ( {\widetilde{h}}^n )_2 (x, y) = 1
\A
D_1 ( {\widetilde{h}}^n )_2 (x, y)
   = n + \sum_{k = 0}^{n - 1} {\widetilde{r}}' (x + k \te),
\]
while
\[
D_1^m ( {\widetilde{h}}^n )_2 (x, y)
   = \sum_{k = 0}^{n - 1} {\widetilde{r}}^{(m)} (x + k \te)
\]
for $m \geq 2$, and all other partial derivatives are zero.
It follows that
\[
\| D_1 ( {\widetilde{h}}^n )_2 \|_{\infty}
  \leq n + n \| {\widetilde{r}}' \|_{\infty}
\A
\| D_1^m ( {\widetilde{h}}^n )_2 \|_{\infty}
  \leq  n \| {\widetilde{r}}^{(m)} \|_{\infty},
\]
while all other partial derivatives are bounded by constants independent
of $n$.
So $h$ is tempered by Lemma~\ref{TCond}.
\end{exa}

We do not know whether the \md s of odd spheres constructed in
\cite{FH} and~\cite{Wn},
the \md\  of $S^2 \times S^1$ of Example~\ref{3DNonHomeo},
or the \md s of $S^n \times S^1$ for odd $n$ of Example~\ref{DiffDims},
can be chosen to be tempered, although it seems reasonable to expect
that they can be~\cite{Kt}.

\begin{qst}\label{IsoImpFlipConj}
Let $M_1$ and $M_2$ be \csmf s,
and let $h_1 \colon M_1 \to M_1$ and $h_2 \colon M_2 \to M_2$
be tempered \md s.
Suppose that the smooth crossed products
${\mathcal{S}} (\Z, \, C^{\infty} (M_1), \, h_1)$
and ${\mathcal{S}} (\Z, \, C^{\infty} (M_2), \, h_2)$ are isomorphic.
Does it follow that $h_1$ and $h_2$ are flip conjugate?
\end{qst}

We have left one point ambiguous in this question: should the
\hme\  implementing the flip conjugacy be required to be smooth?
This makes a difference, even on $S^1$. 
See Theorem 12.5.1 and the preceding discussion in~\cite{KH},
and for further examples also Theorem 12.6.1 in~\cite{KH}.
(We do not know if this kind of behavior can occur for tempered \md s.)

Even if the flip conjugacy is merely required to be continuous,
we suspect that in general isomorphism does not imply flip conjugacy.
On the other hand, isomorphism of the smooth crossed products
is probably a much more restrictive condition
than isomorphism of the \tgca s.
As specific evidence that this might be the case, we offer
Example~\ref{SmNCTorus}.
One way to extract extra information is via the computation of
cyclic cohomology.
This is how the nonisomorphism in Example~\ref{SmNCTorus} was proved
in~\cite{BCEN}.
Cyclic cohomology for crossed products by $\Z$
has been studied in~\cite{Ns}.

\section{Ranks and Schweitzer's example}

\indent
Our consideration of Question~\ref{IsoImpFlipConj}
made us realize just how little
is known about smooth crossed products by tempered \dfe s.
It is known that the smooth crossed product is spectrally invariant
in the \tgca, by Corollary~7.16 of~\cite{Sw1}.
For the significance of this, see the introduction to~\cite{Sw1}
and Section~1 of~\cite{Sw0}.
In particular,
${\mathcal{S}} (\Z, \, C^{\infty} (M), \, h)$ is closed under
holomorphic functional calculus evaluated in $C^* (\Z, M, h)$,
and the inclusion induces an isomorphism on K-theory.
(See~\cite{PS} for more on
the K-theory of smooth crossed products by $\Z$.)
These seem to be the basic properties wanted for noncommutative
differential geometry as in~\cite{Cn}.
Indeed, if it is in fact true that isomorphism of smooth crossed
products is much less common than isomorphism of the C*~crossed
products, then the examples in Section~\ref{Exs} provide examples
of C*-algebras which have quite different natural smooth structures.

However, as far as we can tell, the following questions all remain open,
even for the smooth irrational rotation algebras.
They are motivated by the importance of the stable and real ranks
as invariants for \ca s, and by the role they play in the
classification of crossed product \ca s.

\begin{qst}\label{SmStRank}
Let $M$ be a \csmf, and let $h \colon M \to M$ be a tempered \md.
Does it follow that the smooth crossed product
${\mathcal{S}} (\Z, \, C^{\infty} (M), \, h)$ has stable rank one,
that is, that the invertible elements are dense?
\end{qst}

It is known that $C^* (\Z, M, h)$ has stable rank one
(\cite{LP2}; Corollary~1.2 of~\cite{LP1b}),
but it is not known whether stable rank one passes to a spectrally
invariant subalgebra,
or even a strongly spectrally invariant subalgebra.
The invertible elements of ${\mathcal{S}} (\Z, \, C^{\infty} (M), \, h)$
are of course dense in the topology of $C^* (\Z, M, h)$, but this is
not what is being asked for.

\begin{qst}\label{SmRRank}
Let $M$ be a \csmf, and let $h \colon M \to M$ be a tempered \md.
Suppose $C^* (\Z, M, h)$ has real rank zero.
Does it follow that the selfadjoint invertible elements in
${\mathcal{S}} (\Z, \, C^{\infty} (M), \, h)$ are dense in the
set of all selfadjoint elements?
\end{qst}

\begin{qst}\label{SmFS}
Let $M$, $h$, and $C^* (\Z, M, h)$ be as in Question~\ref{SmRRank}.
Does it follow that the selfadjoint elements in
${\mathcal{S}} (\Z, \, C^{\infty} (M), \, h)$ with finite spectrum
are dense in the set of all selfadjoint elements?
\end{qst}

For \ca s, the properties asked for in Question~\ref{SmRRank}
and Question~\ref{SmFS} are equivalent---both are real rank zero.
We give a dense selfadjoint subalgebra
of $C_0 (\N)$ which is a Banach algebra in its own topology, 
is strongly spectrally invariant in $C_0 (\N)^+$ (even satisfying
the Blackadar-Cuntz differential seminorm conditions~\cite{BC}),
and is closed under $C^{\infty}$ functional calculus
(even $C^1$ functional calculus) in $C_0 (\N)^+$ for
selfadjoint elements, but for which conclusion in Question~\ref{SmFS}
does not hold.
The conclusion in Question~\ref{SmRRank} does hold, so this example
shows two things: that these conditions are not equivalent for
Banach *-algebras, and that the conclusion in Question~\ref{SmFS}
does not pass to strongly spectrally invariant subalgebras.
This example was constructed by Larry Schweitzer, and is reproduced
here with his permission.

\begin{exa}[Schweitzer]\label{SAInvDNotFS}
By convention, we take $\N = \{ 1, 2, \dots \}$.
Set
\[
B = C_0 (\N)^+ = \{ a \in l^{\infty} ( \N ) \colon
    {\mbox{$\limi{n} a (n)$ exists}} \}.
\]
Write $\| \cdot \|_{\infty}$ for its norm.
Let $\ld \colon B \to \C$ be evaluation at $\infty$, that is,
$\ld (a) = \limi{n} a (n)$.
Define
\[
B_0 = \{ a \in B \colon
    {\mbox{$\limi{n} n [a (n) - \ld (a) ]$ exists}} \},
\]
and define $\om \colon B_0 \to \C$ by
$\om (a) = \limi{n} n [a (n) - \ld (a) ]$.
Then for $a \in B_0$ define
\[
\| a \|_{\om} = \sup_{n \in \N} n | a (n) - \ld (a) |
\A \| a \| = \| a \|_{\infty} + \| a \|_{\om}.
\]
\end{exa}

We establish the properties of this example in a sequence of lemmas.

\begin{lem}\label{BasicEst}
Let $a, \, b \in B_0$.
Then
$\| a b \|_{\om}
 \leq \| a \|_{\infty} \| b \|_{\om} + \| a \|_{\om} \| b \|_{\infty}$.
\end{lem}

\begin{pff}
It is obvious that $\ld (a b) = \ld (a) \ld (b)$ and
$| \ld (a) | \leq \| a \|_{\infty}$.
Therefore
\begin{align*}
\| a b \|_{\om} & = \sup_{n \in \N} n | a (n) b (n) - \ld (a) \ld (b) |
         \\
& \leq \sup_{n \in \N} n | a (n) - \ld (a) | \cdot | b (n) |
        + \sup_{n \in \N} n | \ld (a) | \cdot | b (n) - \ld (b) |  \\
& \leq  \| a \|_{\om} \| b \|_{\infty} + \| a \|_{\infty} \| b \|_{\om}.
\end{align*}
This is the result.
\end{pff}

\begin{cor}\label{DSN}
The definitions $T_0 (a) = \| a \|_{\infty}$, $T_1 (a) = \| a \|_{\om}$,
and $T_n (a) = 0$ for $n \geq 2$, give a differential seminorm on
$B_0$ in the sense of Definition~3.1 of~\cite{BC}.
\end{cor}

\begin{pff}
The required inequalities are
\[
\| a b \|_{\infty} \leq \| a \|_{\infty} \| b \|_{\infty} \A
\| a b \|_{\om}
 \leq \| a \|_{\infty} \| b \|_{\om} + \| a \|_{\om} \| b \|_{\infty}.
\]
The first is known and the second is Lemma~\ref{BasicEst}.
\end{pff}

\begin{prp}\label{SubMult}
The norm $\| \cdot \|$ is submultiplicative on $B_0$, satisfies
$\| a^* \|= \| a\|$ for all $a \in B_0$, and satisfies $\| 1 \| = 1$.
\end{prp}

\begin{pff}
For the first part, we estimate:
\begin{align*}
\| a b \| & = \| a b \|_{\infty} + \| a b \|_{\om}
    \leq \| a \|_{\infty} \| b \|_{\infty}
      + \| a \|_{\infty} \| b \|_{\om} + \| a \|_{\om} \| b \|_{\infty}
               \\
 &  \leq \left( \| a \|_{\infty} + \| a \|_{\om} \right)
            \left( \| b \|_{\infty} + \| b \|_{\om} \right).
\end{align*}
The other two parts are obvious.
\end{pff}

\begin{lem}\label{Complete}
The algebra $B_0$ is a Banach *-algebra in $\| \cdot \|$.
\end{lem}

\begin{pff}
It only remains to prove that $B_0$ is complete.
Let $(a_k)$ be a Cauchy sequence in $B_0$.
Then, using $\| a \|_{\infty} \leq \| a \|$ and
$| \om (a) | \leq \| a \|$,
there are $a \in B$ and $\af \in \C$ such that
\[
\limi{k} \| a_k - a \|_{\infty} = 0 \A \limi{k} | \om (a_k) - \af | = 0.
\]
We first show that $a \in B_0$,
by proving that $\limi{n} n [a (n) - \ld (a) ] = \af$.

Let $\ep > 0$.
Choose $N_0$ so large that if $k, \, l \geq N_0$ then
$\| a_k - a_l \| < {\textstyle{\frac{1}{3}}} \ep$.
For such $k$ and $l$, we have in particular
\[
\sup_{n \in \N} n | a_k (n) - a_l (n) - \ld (a_k - a_l) |
  < {\textstyle{\frac{1}{3}}} \ep.
\]
Letting $l \to \infty$, we get
\[
\sup_{n \in \N} n | a_k (n) - a (n) - \ld (a_k - a) |
  \leq {\textstyle{\frac{1}{3}}} \ep
\]
for all $k \geq N_0$.
Now choose $k \geq N_0$ and also so large that
$| \om (a_k) - \af | < {\textstyle{\frac{1}{3}}} \ep$.
For this $k$, choose $N$ so large that if $n \geq N$ then
\[
| n [a_k (n) - \ld (a_k)] - \om (a_k) | < {\textstyle{\frac{1}{3}}} \ep.
\]
For all $n \geq N$ we then have
\begin{align*}
& | n [a (n) - \ld (a)] - \af |  \\
& \hspace*{3em} \mbox{}
   \leq n | a_k (n) - a (n) - \ld (a_k - a) |
     + | n [a_k (n) - \ld (a_k)] - \om (a_k) | + | \om (a_k) - \af | \\
& \hspace*{3em} \mbox{}
   < {\textstyle{\frac{1}{3}}} \ep + {\textstyle{\frac{1}{3}}} \ep
       + {\textstyle{\frac{1}{3}}} \ep = \ep.
\end{align*}

Now we prove that $\| a_k - a \| \to 0$.
Let $\ep > 0$, and choose $N_0$ as before.
As there, for $k \geq N_0$ we get
\[
\sup_{n \in \N} n | a_k (n) - a (n) - \ld (a_k - a) |
  \leq {\textstyle{\frac{1}{3}}} \ep
\A \sup_{n \in \N} | a_k (n) - a (n) |
  \leq {\textstyle{\frac{1}{3}}} \ep,
\]
whence $\| a_k - a \| \leq {\textstyle{\frac{2}{3}}} \ep < \ep$.
\end{pff}

\begin{cor}\label{HolFcnlCalc}
The algebra $B_0$ is closed under holomorphic functional
calculus in $B$.
That is, if $a \in B_0$ and $f$ is a holomorphic function defined
on a neighborhood of $\spec_B (a)$, then $f (a) \in B_0$.
\end{cor}

\begin{pff}
Use Theorem~1.17 of~\cite{Sw1} and Lemma~1.2 of~\cite{Sw0}.
\end{pff}

\begin{prp}\label{CIFcnlCalc}
The algebra $B_0$ is closed under $C^1$ functional
calculus for selfadjoint elements in $B$.
That is, if $a \in B_0$ satisfies $a^* = a$ and $f$ is a $C^1$
function defined on a \nbhd\  of $\spec_B (a)$, then $f (a) \in B_0$.
\end{prp}

\begin{pff}
The element $f (a) \in B$ is given by $f (a) (n) = f (a (n))$ for
all $n \in \N$.
We must prove that this element is in $B_0$, which we do by showing
that
\[
\limi{n} n [f (a (n)) - \ld ( f (a)) ] = f' ( \ld (a)) \om (a).
\]

Let $\ep > 0$.
Choose $\dt > 0$ such that whenever $| t -  \ld (a) | < \dt$ then
\[
| f (t) - f ( \ld (a)) - f' ( \ld (a)) [t - \ld (a)] |
   \leq \left( \frac{\ep}{2 | \om (a) | + 2} \right) | t -  \ld (a) |.
\]
Choose $N \in \N$ such that
$n \geq N$ implies
\[
| a (n) - \ld (a) | < \dt
\A | n [ a(n) - \ld (a)] - \om (a) |
   < \min \left( 1, \, \, \frac{\ep}{2 |f' ( \ld (a)) | + 2 }  \right).
\]
For such $n$, we then use $\ld (f (a)) = f (\ld (a))$ to get
\begin{align*}
& | n [f (a (n)) - \ld ( f (a)) ] - f' (\ld (a) ) \om (a) | \\
& \hspace*{2em} \mbox{}
   \leq \left( \frac{\ep}{2 | \om (a) | + 2} \right)
          \cdot n | a (n) -  \ld (a) | 
     + | f' ( \ld (a)) n [a (n) - \ld (a)] - f' (\ld (a) ) \om (a) | \\
& \hspace*{2em} \mbox{}
   \leq \left( \frac{\ep}{2 | \om (a) | + 2} \right) (|\om (a)| + 1)
     + | f' ( \ld (a)) | \cdot | n [ a (n) - \ld (a)] - \om (a) |  \\
&  < {\textstyle{\frac{1}{2}}} \ep + {\textstyle{\frac{1}{2}}} \ep
   = \ep.
\end{align*}
\end{pff}

The algebra $B_0$ is not closed under \ct\  functional
calculus for selfadjoint elements in $B$.
Take $f (x) = \sqrt{x}$ for $x \geq 0$.
Define $a \in B_0$ by $a (n) = \frac{1}{n}$.
Then $a \in B_0$ but $f (a) \not\in B_0$.

\begin{prp}\label{WeakRank}
The invertible elements of $B_0$ are dense in $B_0$,
and the invertible selfadjoint elements in $B_0$ are dense in the set of
all selfadjoint elements in $B_0$.
\end{prp}

\begin{pff}
It follows from Corollary~\ref{HolFcnlCalc} that every element of
$B_0$ which is invertible in $B$ is also invertible in $B_0$.
So let $a \in B_0$ and let $\ep > 0$.
Choose any real number
$\af \not\in \{ \ld (a) \} \cup \{ a (n) \colon n \in \N \}$
with $| \af | < \ep$.
Then $b = a - \af \cdot 1$ is invertible in $B$ and satisfies
$\| b - a \| < \ep$.
Moreover, if $a$ is selfadjoint, then so is $b$.
\end{pff}

\begin{prp}\label{NoStrongRank}
The selfadjoint elements in $B_0$ which have finite spectrum are not
dense in the set of all selfadjoint elements in $B_0$.
\end{prp}

\begin{pff}
Define $a \in B_0$ by $a (n) = \frac{1}{n}$.
Let $b \in B_0$ have finite spectrum.
Then the range of $b$ is finite,
whence $b (n) = \ld (b)$ for all sufficiently large $n$.
Therefore $\om (b) = 0$.
Since $| \om (b - a) | \leq \| b - a \|_{\om}$ and $\om (a) = 1$,
it follows that $\| b - a \| \geq 1$.
\end{pff}

\end{document}